%
\output={\if N\header\headline={\hfill}\fi
\plainoutput\global\let\header=Y}
\magnification\magstep1
\tolerance = 500
\hsize=14.4true cm
\vsize=22.5true cm
\parindent=6true mm\overfullrule=2pt
\newcount\kapnum \kapnum=0
\newcount\parnum \parnum=0
\newcount\procnum \procnum=0
\newcount\nicknum \nicknum=1
\font\ninett=cmtt9

\font\ninebf=cmbx9

\font\sixbf=cmbx6
\font\ninesl=cmsl9

\font\nineit=cmti9

\font\ninerm=cmr9

\font\sixrm=cmr6
\font\ninei=cmmi9
\font\eighti=cmmi8
\font\sixi=cmmi6
\skewchar\ninei='177 \skewchar\eighti='177 \skewchar\sixi='177
\font\ninesy=cmsy9
\font\eightsy=cmsy8
\font\sixsy=cmsy6
\skewchar\ninesy='60 \skewchar\eightsy='60 \skewchar\sixsy='60
\font\titelfont=cmr10 scaled 1440
\font\paragratit=cmbx10 scaled 1200

\font\name=cmcsc10
\font\emph=cmbxti10

\font\tenmsbm=msbm10
\font\sevenmsbm=msbm7
%

%

%
\font\teneufm=eufm10
\font\seveneufm=eufm7
\font\fiveeufm=eufm5
\newfam\eufmfam
\textfont\eufmfam=\teneufm
\scriptfont\eufmfam=\seveneufm
\scriptscriptfont\eufmfam=\fiveeufm

\font\tenmsam=msam10
\font\sevenmsam=msam7
\font\fivemsam=msam5
\newfam\msamfam
\textfont\msamfam=\tenmsam
\scriptfont\msamfam=\sevenmsam
\scriptscriptfont\msamfam=\fivemsam
\font\tenmsbm=msbm10
\font\sevenmsbm=msbm7
\font\fivemsbm=msbm5
\newfam\msbmfam
\textfont\msbmfam=\tenmsbm
\scriptfont\msbmfam=\sevenmsbm
\scriptscriptfont\msbmfam=\fivemsbm
\def\diag{{\rm diag}}

\def\Ell{{\rm (Ell2)}}
\def\bfE{{\bf E}}
\def\Bbb#1{{\fam\msbmfam\relax#1}}
\def\cz{{\kern0.4pt\Bbb C\kern0.7pt}
}
\def\ez{{\kern0.4pt\Bbb E\kern0.7pt}
}
\def\fz{{\kern0.4pt\Bbb F\kern0.3pt}}
\def\gz{{\kern0.4pt\Bbb Z\kern0.7pt}}
\def\hz{{\kern0.4pt\Bbb H\kern0.7pt}
}
\def\kz{{\kern0.4pt\Bbb K\kern0.7pt}
}
\def\nz{{\kern0.4pt\Bbb N\kern0.7pt}
}
\def\oz{{\kern0.4pt\Bbb O\kern0.7pt}
}
\def\rz{{\kern0.4pt\Bbb R\kern0.7pt}
}
\def\sz{{\kern0.4pt\Bbb S\kern0.7pt}
}
\def\pz{{\kern0.4pt\Bbb P\kern0.7pt}
}
\def\qz{{\kern0.4pt\Bbb Q\kern0.7pt}
}
\newskip\ttglue
\def\ninepoint{\def\rm{\fam0\ninerm}%
  \textfont0=\ninerm \scriptfont0=\sixrm \scriptscriptfont0=\fiverm
  \textfont1=\ninei \scriptfont1=\sixi \scriptscriptfont1=\fivei
  \textfont2=\ninesy \scriptfont2=\sixsy \scriptscriptfont2=\fivesy
  \textfont3=\tenex \scriptfont3=\tenex \scriptscriptfont3=\tenex
  \def\it{\fam\itfam\nineit}%
  \textfont\itfam=\nineit
  \def\sl{\fam\slfam\ninesl}%
  \textfont\slfam=\ninesl
  \def\bf{\fam\bffam\ninebf}%
  \textfont\bffam=\ninebf \scriptfont\bffam=\sixbf
   \scriptscriptfont\bffam=\fivebf
  \def\tt{\fam\ttfam\ninett}%
  \textfont\ttfam=\ninett
  \tt \ttglue=.5em plus.25em minus.15em
  \normalbaselineskip=11pt
  \font\name=cmcsc9
  \let\sc=\sevenrm
  \let\big=\ninebig
  \setbox\strutbox=\hbox{\vrule height8pt depth3pt width0pt}%
  \normalbaselines\rm
  \def\sl{\it}}

\headline={\ifodd\pageno\rightheadline\else\leftheadline\fi}
\def\rightheadline{\ninepoint Paragraphen"uberschrift\hfill\folio}
\def\leftheadline{\ninepoint\folio\hfill Chapter"uberschrift}
\let\header=Y
\def\titel#1{\need 9cm \vskip 2truecm
\parnum=0\global\advance \kapnum by 1
{\baselineskip=16pt\lineskip=16pt\rightskip0pt
plus4em\spaceskip.3333em\xspaceskip.5em\pretolerance=10000\noindent
\titelfont Chapter \uppercase\expandafter{\romannumeral\kapnum}.
#1\vskip2true cm}\def\leftheadline{\ninepoint
\folio\hfill Chapter \uppercase\expandafter{\romannumeral\kapnum}.
#1}\let\header=N
}
\def\Titel#1{\need 9cm \vskip 2truecm
\global\advance \kapnum by 1
{\baselineskip=16pt\lineskip=16pt\rightskip0pt
plus4em\spaceskip.3333em\xspaceskip.5em\pretolerance=10000\noindent
\titelfont\uppercase\expandafter{\romannumeral\kapnum}.
#1\vskip2true cm}\def\leftheadline{\ninepoint
\folio\hfill\uppercase\expandafter{\romannumeral\kapnum}.
#1}\let\header=N
}
\def\need#1cm {\par\dimen0=\pagetotal\ifdim\dimen0<\vsize
\global\advance\dimen0by#1 true cm
\ifdim\dimen0>\vsize\vfil\eject\noindent\fi\fi}
\def\neupara#1{\par\penalty-2000
\procnum=0\global\advance\parnum by 1
\vskip1cm\noindent{\paragratit \the\parnum. #1}%
\def\rightheadline{\ninepoint\S\the\parnum.\ #1\hfill \folio}%
\vskip 8mm\noindent}
\def\Proclaim #1 #2\finishproclaim {\bigbreak\noindent
{\bf#1\unskip{}. }{\it#2}\medbreak\noindent}
%
\gdef\proclaim #1 #2 #3\finishproclaim {\bigbreak\noindent%
\global\advance\procnum by 1
{%
{\relax\ifodd \nicknum
\hbox to 0pt{\vrule depth 0pt height0pt width\hsize
   \quad \ninett#3\hss}\else {}\fi}%
\bf\the\parnum.\the\procnum\ #1\unskip{}. }
{\it#2}
\immediate\write\num{\string\def
 \expandafter\string\csname#3\endcsname
 {\the\parnum.\the\procnum}}
\medbreak\noindent}
\newcount\stunde \newcount\minute \newcount\hilfsvar
\def\uhrzeit{
    \stunde=\the\time \divide \stunde by 60
    \minute=\the\time
    \hilfsvar=\stunde \multiply \hilfsvar by 60
    \advance \minute by -\hilfsvar
    \ifnum\the\stunde<10
    \ifnum\the\minute<10
    0\the\stunde:0\the\minute~Uhr
    \else
    0\the\stunde:\the\minute~Uhr
    \fi
    \else
    \ifnum\the\minute<10
    \the\stunde:0\the\minute~Uhr
    \else
    \the\stunde:\the\minute~Uhr
    \fi
    \fi
    }
\def\Sch{{\rm Sch}}

 \def\calB{{\cal B}}

\def\calM{{\cal M}} 
\def\calO{{\cal O}} \def\calP{{\cal P}}

\def\GL{\mathop{\rm GL}\nolimits}

\def\Hom{\mathop{\rm Hom}\nolimits}

\def\kernel{\mathop{\rm kernel}\nolimits}

\def\mod{\mathop{\rm mod}\nolimits}

\def\Pic{\mathop{\rm Pic}\nolimits}

\def\proj{\mathop{\rm proj}\nolimits}

\def\SL{\mathop{\rm SL}\nolimits}

\def\Spec{\mathop{\rm Spec}\nolimits}

\def\boxit#1{
  \vbox{\hrule\hbox{\vrule\kern6pt
  \vbox{\kern8pt#1\kern8pt}\kern6pt\vrule}\hrule}}
\def\Boxit#1{
  \vbox{\hrule\hbox{\vrule\kern2pt
  \vbox{\kern2pt#1\kern2pt}\kern2pt\vrule}\hrule}}

\def\smallni{\smallskip\noindent }
\def\medni{\medskip\noindent }

\def\Isom{\mathop{\;{\buildrel \sim\over\longrightarrow }\;}}
\def\lo{\longrightarrow}

\def\loma{\longmapsto}
\def\betr#1{\vert#1\vert}

\def\imag{{\rm i}}
\def\pii{\pi {\rm i}}

\def\square{\hbox{\hbox to 0pt{$\sqcup$\hss}\hbox{$\sqcap$}}}
\def\qed{\ifmmode\square\else{\unskip\nobreak\hfil
\penalty50\hskip3em\null\nobreak\hfil\square
\parfillskip=0pt\finalhyphendemerits=0\endgraf}\fi}
\def\pn{\the\parnum.\the\procnum}
\def\downmapsto{{\buildrel
        {\vbox{\hbox{\hskip.2pt$\scriptstyle-$}}}
        \over{\raise7pt\vbox{\vskip-4pt\hbox{$\textstyle\downarrow$}}}}}
\def\DSP{\displaystyle}
\def\PET{\mathop\mid}                
\def\PETER_#1{{\DSP \PET_{#1}}}      
\input box3.num
\nopagenumbers
\immediate\newwrite\num
\nicknum=0  
\let\header=N
\def\transpose#1{\kern1pt{^t\kern-1pt#1}}%

\immediate\openout\num=box3.num
\immediate\newwrite\num\immediate\openout\num=box3.num
\def\RAND#1{\vskip0pt\hbox to 0mm{\hss\vtop to 0pt{%
  \raggedright\ninepoint\parindent=0pt%
  \baselineskip=1pt\hsize=2cm #1\vss}}\noindent}
\noindent
\centerline{\titelfont Parametrization of the box variety by theta functions}%
\def\leftheadline{\ninepoint\folio\hfill
The box variety}%
\def\rightheadline{\ninepoint Introduction\hfill \folio}%
\headline={\ifodd\pageno\rightheadline\else\leftheadline\fi}
\vskip1.5cm
\leftline{\it \hbox to 6cm{Eberhard Freitag\hss}
Riccardo Salvati
Manni  }
  \leftline {\it  \hbox to 6cm{Mathematisches Institut\hss}
Dipartimento di Matematica, }
\leftline {\it  \hbox to 6cm{Im Neuenheimer Feld 288\hss}
Piazzale Aldo Moro, 2}
\leftline {\it  \hbox to 6cm{D69120 Heidelberg\hss}
 I-00185 Roma, Italy. }
\leftline {\tt \hbox to 6cm{freitag@mathi.uni-heidelberg.de\hss}
salvati@mat.uniroma1.it}
\vskip1cm
\centerline{\paragratit \rm  2013}%
\vskip5mm\noindent%
\let\header=N%
\def\imag{{\rm i}}%
{\paragratit Introduction}%
\medni
We consider the graded algebra (the generators have weight one)
$$B=\qz[Z_1,Z_2,Z_3,W_1,W_2,W_3,C]$$
with defining relations
$$\eqalign{W_1^2+W_2^2&=Z_3^2,\cr
W_1^2+W_3^2&=Z_2^2,\cr
W_2^2+W_3^2&=Z_1^2,\cr
W_1^2+W_2^2+W_3^2&=C^2.\cr}$$
This is a normal graded algebra. The associated projective
variety $\proj(B)$
is called the box variety. It is absolutely irreducible.
We denote its complexification by
$$\calB:=\proj(B\otimes_\qz\cz).$$
It is a surface that characterizes
cuboids. The variables $W_i$ give the edges of the cuboid, the variables
$Z_i$ the diagonals of the faces and $C$ the long diagonal.
We mention that there is an unsolved problem, raised by Euler, whether the
box variety contains non-trivial rational points or not. For more details on the
box variety we refer to [vL] and [ST].
\smallskip
In this note we describe a parametrization of the 
box variety  (variety of cuboids) by theta functions. This will imply
that the box variety is a quotient of the product 
$\overline{\hz/\Gamma[8]}\times\overline{\hz/\Gamma[8]}$ of two modular curves
of level 8 by a group of order 8 which comes from
the diagonal action of $\Gamma[4]$. Actually this parametrization
can be defined over
the Gauss number field $K=\qz(\imag)$. We found this parametrization
from an observation of D.~Testa that the box variety can be embedded into a
certain Siegel modular variety which has been described by van Geemen and
Nygaard. This background is not necessary for our note and we will not describe
it here. But we want to point out that this still unpublished work of Testa is
behind the scenes and we are very grateful that Testa explained to us
details of this work.
\smallskip
This parametrization can be used to derive quickly known properties and also
some new ones of the box variety. For example, we give in Sect.~2 a modular description of
the automorphism group. It can be realized through a  
subgroup of $\SL(2,\gz)\times\SL(2,\gz)$. In [ST], [vL] 140 rational and elliptic curves
of the minimal model of the box variety which give generators of the
Picard group have been described. We  describe them
in Sect.~3 in a very simple way as certain modular curves.
\smallskip
In Sect.~3 we consider smooth curves in the box variety. We prove an estimate that
shows how their genus grows with their degree. As a consequence, smooth
rational and elliptic curves have a bounded degree. This can be considered as
validation of a conjecture made in [ST] that the 140  curves
described in [vL] exhaust all  rational and elliptic
curves. This has been proved in [ST] for degrees
$\le 4$.
\smallskip
In the paper [Be] of 
Beauville the box variety arises as a member of a whole family having the same
properties, namely to be complete intersections of 4 quadrics in $\pz^6$ with an even set
of 48 nodes. In this paper Beauville also
describes a certain smooth two-fold Galois covering $X$ of the box variety. 
It is unramified outside the 48 nodes and it is a minimal surface of general type
with $q=4$, $p_g=7$, $K^2=32$.
In Sect.~4 we give
a very simple modular description of it.
\smallskip
In Sect.~5 a certain involution $\sigma$ of the box variety $\calB$ is considered.
We use the modular description to realize the quotient $\calB/\sigma$ 
as a Kummer variety.
\smallskip
In the last section we consider a certain moduli problem which gives the
realization of the box variety as fine moduli scheme over $\qz(\imag)$
classifying pairs $(E,F)$ of elliptic curves with level 4 structures and a
compatible  isomorphism
$E[8]\to F[8]$. This is closely related to work of E.~Kani [Ka].
\smallskip
We want to thank A.~Beauville, E.~Kani and D.~Testa for helpful discussions.
\neupara{Generalities about modular groups}%
We use the standard notations 
$$\Gamma[N]=\kernel(\SL(2,\gz)\lo\SL(2,\gz/N\gz))$$
for the principal congruence
subgroup of level $N$ of the elliptic modular group and 
$$\eqalign{\Gamma_0[N]=\big\{\,M=\pmatrix{a&b\cr c & d}\in
\SL(2,\gz);\ c\equiv 0\mod N\,\big\},\cr
\Gamma_1[N]=\big\{\,M=\pmatrix{a&b\cr c & d}\in
\SL(2,\gz);\quad a\equiv b\equiv 1\mod N,\ c\equiv 0\mod N\,\big\}.\cr}$$
We also will use the Igusa groups
$$\Gamma[N,2N]=\big\{\, M=\pmatrix{a&b\cr c&d};\quad ab\equiv cd\equiv 0\;\mod\;2N\,\big\}.$$
In the following we define $\sqrt a$ for a non-zero complex number by the principal
part of the logarithm. This means that the real part is positive if $a$ is not real and negative
and that $\sqrt a=\imag\>\strut_{\hbox{\sevenrm +}}\hskip-1mm\sqrt{\betr a}$ if $a$ is real and negative. 
Let $\Gamma\subset\SL(2,\gz)$ be a subgroup of finite index
and let $r$ be an integral 
number. By a {\it multiplier system\/} of weight $r/2$ one understands a map $v:\Gamma\to S^1$
such that 
$$v(M)\sqrt{c\tau+d}^{\,r},\qquad M=\pmatrix{a&b\cr c&d},$$ 
is a cocycle. Then the space
$[\Gamma,r/2,v]$ of entire modular forms can be defined in the usual way. Their
transformation law is 
$$f(M\tau)=v(M)\sqrt{c\tau+d}^{\,r}f(\tau).$$
There are two basic multiplier systems. The theta multiplier system $v_\vartheta$ is a
multiplier system of weight $1/2$ on the theta group 
$$\Gamma_\vartheta:=\Gamma[1,2].$$ 
It can be defined as the  multiplier system
of the theta function
$$\vartheta(\tau)=\sum_{m=-\infty}^\infty e^{\pii n^2\tau}.$$
The theta group is generated by ${1\,2\choose0\,1}$ and ${0\,-1\choose 1\;\;0}$.
From the  theta inversion formula 
$\vartheta(-1/\tau)=\sqrt{\tau/\imag}\vartheta(\tau)$ we get
$$v_\vartheta\pmatrix{1&2\cr0&1}=1,\quad v_\vartheta\pmatrix{0&-1\cr1&0}=e^{-\pii /4}.$$
We also have to consider the theta function of 
second kind $\Theta(\tau):=\vartheta(2\tau)$.
This is a modular form for $\Gamma_0[4]$. We denote its multiplier system by $v_\Theta$.
Both multiplier systems $v_\vartheta,v_\Theta$ agree on $\Gamma[8]$.
\smallskip
For given $\Gamma$, $r_0$, $v$ we can consider the graded algebra
$$A(\Gamma,r_0,v):=\sum_{r\in\gz}[\Gamma,rr_0,v^r].$$
If it is clear which $(r_0,v)$ has to be considered we will simply write
$A(\Gamma)$ for this algebra. This is a finitely generated algebra of Krull dimension 2.
The associated projective curve $\proj(A(\Gamma))$ can be identified with
$$\overline{\hz/\Gamma}=\hz^*/\Gamma\quad\hbox{where}\quad
\hz^*=\hz\cup\qz\cup\{\infty\}.$$
We have to consider more generally subgroups of finite index $\Gamma\subset
\SL(2,\gz)\times\SL(2,\gz)$. A multiplier system $v$ of weight $r/2$ now means a function
$v:\Gamma\to S^1$ such that 
$$v(M_1,M_2)\sqrt{c_1\tau+d_1}^{\,r}\sqrt{c_2\tau+d_2}^{\,r}$$
is a cocycle. The spaces of modular forms $[\Gamma,r/2,v]$ (now functions of
two variables) and the algebras $A(\Gamma)=A(\Gamma,r_0,v)$ are defined in the obvious way.
\smallskip
Let $N$ be a divisor of the natural number $N'$. In this paper the group
$$\Delta(N,N')=\big\{\,(M_1,M_2)\in\Gamma[N]\times\Gamma[N],
\quad M_1\equiv M_2\;\mod\; N'\,\big\}$$
will play a role. It is generated by $\Gamma[N']\times\Gamma[N']$ and the diagonally
embedded $\Gamma[N]$.
\neupara{A parametrization of the box variety by theta functions}%
We make use of the Jacobi theta functions
$$\vartheta_{a,b}(z)=\sum_{n=-\infty}^\infty e^{\pii(n+a/2)^2z+b(n+a/2)}.$$
Here $(a,b)$  is one of the three pairs $(0,0)$, $(1,0)$, $(0,1)$.
These three functions are modular forms of weight $1/2$ with respect to the three conjugate
groups of the theta group. The multiplier systems agree on the group $\Gamma[4,8]$ with
$v_\vartheta$. We also will consider the two theta functions of the second kind
$\vartheta_{00}(2\tau)$ and $\vartheta_{10}(2\tau)$. They are modular forms for the two conjugated
groups of $\Gamma_0[4]$. Their multiplier systems agree on $\Gamma[2,4]$ with $v_\Theta$.
\smallskip
So we see that the 5 functions
$$\vartheta_{00}(z),\quad\vartheta_{10}(z),\quad
\vartheta_{01}(z),\quad\vartheta_{00}(2z),\quad
\vartheta_{10}(2z)$$
have the same multiplier system on $\Gamma[8]$.
Hence they are contained in the ring
$$A(\Gamma[8]):=\bigoplus_{r\in\gz}[\Gamma[8],r/2,v_\vartheta^r].$$
It is not difficult to show the following result.
The details have been worked out in the Heidelberg Diplomarbeit [Br].
\proclaim
{Theorem}
{One has
$$A(\Gamma[8])=\cz[\vartheta_{00}(z),\vartheta_{10}(z),
\vartheta_{01}(z),\vartheta_{00}(2z),
\vartheta_{10}(2z)].$$
Defining relations are the classical theta relations
$$\eqalign{
\vartheta_{00}(z)^2&=\vartheta_{00}(2z)^2+\vartheta_{10}(2z)^2,\cr
\vartheta_{01}(z)^2&=\vartheta_{00}(2z)^2-\vartheta_{10}(2z)^2,\cr
\vartheta_{10}(z)^2&=2\vartheta_{00}(2z)\vartheta_{10}(2z).\cr}$$
}
StrE%
\finishproclaim
Since the multiplier system $v_\vartheta$ is defined on the theta group $\Gamma_\vartheta$,
we can define an action of the theta group on $A(\Gamma[8])$ by the formula
$$f\vert M(\tau)= v_\vartheta(M)^{-r} \sqrt{cz+d}^{-r}f(Mz).$$
This is an action from the right, $f\vert (M_1M_2)=(f\vert M_1)\vert M_2$.
We describe it by means of matrices. For this we have to use the action of $\GL(n,\cz)$
on a complex vector space $V$ with basis $e_1,\dots,e_n$ from the right. It is defined by
$Ae_i=\sum a_{ij}e_j$. If we write an element of $V$ in the form $\sum x_ie_i$ then  this
means that the row $x=(x_1,\dots,x_n)$ had to be multiplied from the right by the matrix $A$.
Standard theta transformation formulas give the following result.
\proclaim
{Lemma}
{The matrix ${1\,2\choose0\,1}$ acts with respect to the basis
$$\vartheta_{00}(z),\quad\vartheta_{10}(z),\quad
\vartheta_{01}(z),\quad\vartheta_{00}(2z),\quad
\vartheta_{10}(2z)$$
through the diagonal matrix with the diagonal entries
$$1,\imag,1,1,-1.$$
The matrix ${0\,-1\choose1\,\phantom{-}0}$ acts with respect to this basis
through the matrix
$$\pmatrix{1&0&0&0&0\cr
0&0&1&0&0\cr
0&1&0&0&0\cr
0&0&0&1/\sqrt2&1/\sqrt2\cr
0&0&0&1/\sqrt2&-1/\sqrt2\cr}$$
}
ActGt%
\finishproclaim
We are interested in the action of $\Gamma[4]$ on $A(\Gamma[8])$.
The factor group $\Gamma[4]/\Gamma[8]$ is ismorphic to $\gz/2\gz^3$. It is generated by
the images of the matrices
$$T=\pmatrix{1&4\cr 0& 1},\quad T'=\pmatrix{1&0\cr 4& 1},\quad
R=\pmatrix{5&8\cr 8& 13}.$$
From Lemma \ActGt\ we get the following result.
\proclaim
{Lemma}
{The generators $T,T',R$ of $\Gamma[4]/\Gamma[8]$ act on $A(\Gamma[4])$ by means
of the diagonal matrices
$$\eqalign{
T&\loma\diag(1,-1,1,1,1),\cr
T'&\loma\diag(1,1,-1,,1,1),\cr
R&\loma\diag(1,1,1,-1,-1).\cr}$$
}
ActGv%
\finishproclaim
Now we consider modular forms of two variables. We consider the ring
$A(\Gamma[8]\times\Gamma[8])$ of modular forms of integral or half integral
weight $r/2$ with respect to the multiplier system $(v_\vartheta(M_1)v_\vartheta(M_2))^r$.
It is clear that
$$A(\Gamma[8]\times\Gamma[8]):=\cz[f(z)g(w)],\quad 
f,g\in\{\vartheta_{00}(\cdot),\vartheta_{10}(\cdot),
\vartheta_{01}(\cdot),\vartheta_{00}(2\cdot),
\vartheta_{10}(2\cdot)\}.$$
We want to determine the subring $A(\Delta(4,8))$ of modular forms
with respect to the group $\Delta(4,8)$. This is the ring of invariants
with respect to the diagonal action of $\Gamma[4]$ by means of the action
$$f(z,w)\loma v_\vartheta(M)^{-2r}\sqrt{cz+d}^{\,-r}\sqrt{cw+d}^{\,-r}f(Mz,Mw).$$
Using Lemma \ActGv\ it is obvious that the forms
$$\eqalign{
&\vartheta_{00}(z)\vartheta_{00}(w),\cr
&\vartheta_{10}(z)\vartheta_{10}(w),\cr
&\vartheta_{01}(z)\vartheta_{01}(w),\cr
\noalign{\vskip1mm}
&\vartheta_{00}(2z)\vartheta_{00}(2w),\cr
&\vartheta_{00}(2z)\vartheta_{10}(2w),\cr
&\vartheta_{10}(2z)\vartheta_{00}(2w),\cr
&\vartheta_{10}(2z)\vartheta_{10}(2w).\cr
}$$
are invariant under
$\Gamma[4]$. Moreover, one can show that they generate the invariant ring.
\proclaim
{Theorem}
{There is an isomorphism
$$B\otimes_\qz \cz\Isom A(\Delta(4,8))$$
which is given by
$$\eqalign{
Z_1&\loma\vartheta_{01}(z)\vartheta_{01}(w),\cr
Z_2&\loma\vartheta_{00}(z)\vartheta_{00}(w),\cr
Z_3&\loma\vartheta_{10}(z)\vartheta_{10}(w),\cr
\noalign{\vskip1mm}
C&\loma \vartheta_{00}(2z)\vartheta_{00}(2w)+
\vartheta_{10}(2z)\vartheta_{10}(2w),\cr
\noalign{\vskip1mm}
W_1&\loma \vartheta_{10}(2z)\vartheta_{00}(2w)+
\vartheta_{00}(2z)\vartheta_{10}(2w),\cr
W_2&\loma \imag(\vartheta_{10}(2z)\vartheta_{00}(2w)-
\vartheta_{00}(2z)\vartheta_{10}(2w)),\cr
W_3&\loma \vartheta_{00}(2z)\vartheta_{00}(2w)-
\vartheta_{10}(2z)\vartheta_{10}(2w).\cr
}$$
Hence we have $\calB\cong \overline{\hz\times\hz/\Delta(4,8)}$.}
MIs%
\finishproclaim
{\it Proof.\/} Classical theta relations
$$\eqalign{
\vartheta_{00}(z)^2&=\vartheta_{00}(2z)^2+\vartheta_{10}(2z)^2,\cr
\vartheta_{01}(z)^2&=\vartheta_{00}(2z)^2-\vartheta_{10}(2z)^2,\cr
\vartheta_{10}(z)^2&=2\vartheta_{00}(2z)\vartheta_{10}(2z).\cr}$$
show that this is
a homomorphism. Obviously it is surjective.
Since $A(\Delta(4,8))$ is an integral domain of Krull
dimension three and since $B$ also has dimension three,
this homomorphism must be an isomorphism.\qed
\smallskip
The modular picture can be used to recover known properties of the box variety.
We mention some of them.
\smallskip
First we describe the automorphism group of the box variety. The group $\Delta(4,8)$
is a normal subgroup of $\Delta(1,2)$. The index is 768. Hence the quotient
$\Delta(1,2)/\Delta(4,8)$ is a subgroup of order 768 of the automorphism group.
The involution $(z,w)\mapsto(w,z)$ gives an extra automorphism. Both together generate
a subgroup of order 1$\,$536 of the automorphism group. Due to [ST] the order
of the automorphism group is 1$\,$536. Hence we described the full automorphism group.
\smallskip
Now we describe the singularities.
It is known that the box variety has 48 singularities which all are nodes. In the modular picture
they correspond to some zero dimensional cusps. 
These are the images of the points $(a,b)$ where $a,b\in\qz\cup\{\infty\}$.
There are two types of such points. It may happen that $(a,b)$ is the fixed point
of pair $(M_1,M_2)$ of parabolic elements. The typical case is
$(\infty,\infty)$ and $A=B={1\,4\choose 0\,1}$. The group $\Delta(1,2)$ 
acts transitively on them. There are pairs which do not have this property.
The precise picture is as follows.
\proclaim
{Proposition}
{The box variety $\overline{\hz\times \hz/\Delta(4,8)}$ contains $96$
zero dimensional cusps. They decompose into two orbits of $48$ cusps under $\Delta(1,2)$.
The orbit containing the image of $(\infty,\infty)$ defines the singular locus.}
SingL%
\finishproclaim
A slightly different way to see this is to consider the Galois coverings
$$\overline{\hz/\Gamma[8]}\times \overline{\hz/\Gamma[8]}\lo 
\overline{\hz\times \hz/\Delta(4,8)}
\lo\overline{\hz/\Gamma[2]}\times \overline{\hz/\Gamma[2]}.$$
The covering group of the first cover is
$G=\Delta(4,8)/\Gamma[8]\times\Gamma[8]\cong(\gz/2\gz)^3$. 
The singular points of
the box variety are the images of the fixed points of $G$. 
They agree with the fibres of three zero dimensional cusps of 
$\overline{\hz/\Gamma[2]}\times \overline{\hz/\Gamma[2]}$ which are of the form
$(a,a)$. Thy can be represented by $(\infty,\infty)$, $(0,0)$ and $(1,1)$.
In the typical case
$(\infty,\infty)$ one can take  $p=e^{2\pii z/8}$, $q=e^{2\pii w/8}$
as uniformizing parameters of 
$\overline{\hz/\Gamma[8]}\times \overline{\hz/\Gamma[8]}$ . The stabilizer
in $G$ is generated by the translation $(z,w)\mapsto (z+4,w+4)$ wich acts by
$(p,q)\mapsto-(p,q)$. Hence the singularity appears as quotient singularity
of the type
$(\cz\times\cz)/\pm$ which actually is a node.
\smallskip
We denote by $\tilde\calB$ the minimal resolution of the 48 nodes. The exceptional divisor
is the union of 48 lines.
\smallskip
Next we describe the holomorphic differential forms on $\tilde\calB$.
The modular curve $\overline{\hz/\Gamma[8]}$ has genus 5. The differentials
$$\eqalign{
\omega_1(z)&=\vartheta_{00}(z)^2\vartheta_{01}(z)\vartheta_{10}(z)dz,\cr
\omega_2(z)&=\vartheta_{00}(z)\vartheta_{01}(z)^2\vartheta_{10}(z)dz,\cr
\omega_3(z)&=\vartheta_{00}(z)\vartheta_{01}(z)\vartheta_{10}(z)^2dz.\cr
\omega_4(z)&=\vartheta_{00}(2z)\vartheta_{00}(z)\vartheta_{01}(z)\vartheta_{10}(z) dz.\cr
\omega_5(z)&=\vartheta_{10}(2z)\vartheta_{00}(z)\vartheta_{01}(z)\vartheta_{10}(z)dz.\cr
}$$
are holomorphic on
$\overline{\hz/\Gamma[8]}$, since the defining modular forms are  cusp forms.
A simple computation gives that
\smallni
$
\psi_1=\omega_1(z)\wedge \omega_1(w),\;\psi_2=\omega_2(z)\wedge \omega_2(w),\;\psi_3=\omega_3(z)\wedge \omega_3(w),\;
\psi_4=\omega_4(z)\wedge \omega_4(w),\;\psi_5=\omega_4(z)\wedge 
\omega_5(w),\;\psi_6=\omega_5(z)\wedge \omega_4(w),\; \psi_7=\omega_5(z)\wedge \omega_5(w)$
\smallni
are  $\Delta(4,8)$-invariant holomorphic differential forms.
One can check that
they extend holomorphically to the desingularization $\tilde\calB$.
In this way one can recover the result of [ST]
that the minimal resolution of the box variety has geometric genus 7.
One can also derive from this picture that the box variety is of general type.
\smallskip
In the paper [ST] the structure of the Picard group of $\tilde\calB$
has been determined.  It is a free 
abelian group of rank 64. Stoll and Testa proved that certain 140 curves defined
already in [vL] generate this group. There are 80 rational and
60 elliptic curves. 
We first describe the rational curves. The 48 exceptional curves belong to them.
The remaining 32 rational curves have the following easy modular description.
\proclaim
{Proposition}
{The  equations
$w=Mz+k$ where $M$ runs through a system of representatives of $\Gamma[4]/\Gamma[8]$
and $k\in\{0,2,4,6\}$ define $32$ smooth rational curves in the box variety.
There union is the zero set of the modular form
$$\vartheta_{00}(z)^4\vartheta_{01}(w)^4-\vartheta_{01}(z)^4\vartheta_{00}(w)^4\qquad
(=4\imag W_1W_2W_3C).$$
}
SmRat%
\finishproclaim
Next we describe the elliptic curves. Part of them is in the Satake boundary.
The Satake boundary is the union of the images of
$\hz^*\cup\{a\}$ and $\{a\}\cup\hz^*$, where $a\in\qz\cup\{\infty\}$.
It is easy to work out the structure.
\proclaim
{Proposition}
{The 
Satake boundary consists of $12$ (smooth) elliptic curves. Each of them
contains $8$ singular and $8$ smooth zero dimensional cusps.
The whole Satake boundary is the zero set of the modular form
$$\vartheta_{00}(z)\vartheta_{10}(z)\vartheta_{01}(z)
\vartheta_{00}(w)\vartheta_{10}(w)\vartheta_{01}(w)\qquad(=Z_1Z_2Z_3).$$
}
EllSat%
\finishproclaim
We will not give the details of the proof but we explain a typical boundary curve.
We take the image of $\hz^*\times\{\infty\}$. 
This is the modular curve with respect to the group
$\Gamma_1[8]\cap\Gamma[4]$. 
It contains $\Gamma[8]$ as a  subgroup of index 2. It is not difficult 
to work out the structure of the ring of modular forms.
Details can be found in [Kl].
\proclaim
{Proposition}
{The ring $A(\Gamma_1[8]\cap\Gamma[4])$ of all modular forms of half 
integral weight for the group $\Gamma_1[8]\cap\Gamma[4]$ with
respect to the multiplier system $v_\vartheta^r$ is generated by
$$a=\vartheta_{0,0}(z),\quad b=\vartheta_{0,1}(z),\quad
c=\vartheta_{0,0}(2z),\quad d=\vartheta_{1,0}(2z).$$
Defining relations are
$$a^2=c^2+d^2,\quad
b^2=c^2-d^2.$$
}
KleinR%
\finishproclaim
This is an intersection of two quadrics in $\pz^3$ and hence an
elliptic curve.
So this describes one of the 12 elliptic curves in the boundary part of the
box variety.
\smallskip
Finally we describe the 48 elliptic curves that are not contained in the boundary.
\proclaim
{Proposition}
{The equation $w=z+1$ describes an elliptic curve in the box variety which 
also can be defined by the equations
$$W_1=W_2,\
Z_1=Z_2,\
\sqrt2W_1=Z_3,\
W_3^2 + Z_3^2 - C^2=0,\
2Z_2^2 + Z_3^2 -2 C^2=0.
$$
Applying the group $\Delta(1,2)$ one gets $48$ elliptic curves.
}
EllFi%
\finishproclaim
\neupara{Curves in the box variety}%
We want to study irreducible curves  $C\subset\calB$ in the box variety. 
We have to make the rather strong assumption that 
the normalization map $\bar C\to C$ is bijective. We denote by $g$ the genus of $\bar C$.
We use the following fact which has been explained in [ST] and [Be] and which can be seen
from the explicit description of the holomorphic 2-forms on $\tilde\calB$ above: 
the canonical map (defined by the canonical divisor on $\tilde\calB$) is the composition of
the natural projection $\tilde\calB\to\calB$ and the original embedding $\calB\to \pz^6$.
This shows the that the degree $d$ of $C$ in $\pz^6$ equals the intersection 
number of the strict transform
of $C$ in $\tilde\calB$ and a canonical divisor on $\tilde\calB$.
\proclaim
{Theorem}
{Let $C\subset\calB$ be a curve such that the normalization map $\bar C\to C$ is bijective.
Let $g$ be the genus of $\bar C$ and $d$ be the degree of $C$. Then the inequality
$$d\le 176+16g$$
holds.}
DegGen%
\finishproclaim
As a consequence, rational end elliptic curves have bounded degree. This supports a conjecture
in [ST] that each rational or elliptic curve in $\tilde \calB$ belongs to the system of 140 
curves described
above.
\smallni
{\it Proof of Theorem \DegGen.\/} 
Let $k$ be a natural number.
We consider a modular form of weight $4k$ for the group $\Delta(4,8)$. Then we
consider the tensor
$$T=\Delta(z)^k\Delta(w)^kf(z,w)(dzdw)^{8k}$$
in the algebra of symmetric tensors. Since the modular form $\Delta$ has weight 12, this tensor is invariant under
$\Delta(4,8)$.
Hence it  defines a meromorphic tensor on $\tilde\calB$. 
Using the coordinates $p=e^{2\pii z/8}$, $q=e^{2\pii w/8}$,
it is easy to check that this
tensor is holomorphic outside the 48 exceptional curves. In the exceptional curves it may have
poles. We can lift the curve $\bar C$ to a holomorphic map $\varphi:\bar C\to\tilde\calB$.
Then we consider the pulled back tensor $\varphi^*T$. This is a meromorphic tensor of degree $16k$
on $\bar C$. 
\smallskip
We can assume that $C$ is not the image of a $\hz^*\times \{a\}$ or $\{a\}\times\hz^*$ for an $a\in\hz^*$,
since for these curves the theorem can be seen directly.
Then the tensor $\varphi^*T$
vanishes if and only $f$ vanishes
along $C$ as a function. Since the weight $k$ can be made large we can choose $f$ that it doesn't vanish
along $C$ and in addition we can get that $f$ does not vanish at any of the 48 nodes in $\calB$.
\smallskip
We assume that $C$ contains one the nodes, for example the image of the cusp $(\infty,\infty)$.
Then we can consider a parametrization of the curve close to this cusp. 
A simple lemma (compare [Fr], Satz 1) shows
that the curve can be parametrized by a holomorphic map $\alpha:\hz\to\hz\times\hz$ as follows: 
$$\alpha(z)=z(a_1,a_2)+(\Phi_1(q),\Phi_2(q)),\quad q=e^{2\pii\tau},$$
where $\Phi_i$ are holomorphic at $q=0$. The pair $(a_1,a_2)$ is contained in the
translation lattice, i.e. 
$$a_1\equiv a_2\equiv 0\;\mod\; 4,\quad a_1+ a_2\equiv 0\;\mod\; 8.$$
The numbers $a_1,a_2$ are both positive. This implies $a_1+a_2\ge 8$.
\smallskip
We study the poles and zeros of the tensor $\varphi^*T$. Since its divisor is the $16k$-multiple of a canonical
divisor, we have
$$16(2g-2)k=\#\hbox{zeros}-\#\hbox{poles}.$$
First we estimate the number of poles of $\varphi^*T$ from above (counted with multiplicity).
The poles are intersection points of $\bar C$ with the exceptional divisor. Since $\bar C\to C$ is bijective,
$\bar C$ can meet each of the 48 exceptional curves at most once. Hence the set of poles contains at
most 48 points. 
\smallskip
We have to estimate the pole order. It is sufficient to do this for
the standard node, i.e.~the image of the cusp $(\infty,\infty)$ and we can do this by means of the curve
lemma.
The term $(dz dw)^{8k}$ contributes with $16k$ to the pole order and
$\Delta(z)\Delta(w)$ contributes with a zero of order $(a_1+a_2)k\ge 8k$. Hence the pole order of the
tensor at the node is at most $8k$. Since we have 48 nodes the total pole order is estimated by
$384k$.
\smallskip
Next we estimate the number of zeros from below. 
Each intersection point of the zero divisor of $f$ with the curve produces
a zero. (Since we assumed that $f$ doesn't vanish at the nodes, there is no conflict with the poles of $T$.)
Since the zero divisor of $f$ is a $2k$-multiple of the canonical divisor, we get that there are
at least $2kd$ zeros. So we get
$$16(2g-2)k=\#\hbox{zeros}-\#\hbox{poles}\ge 2kd-384 k.$$
This finishes the proof of Theorem \DegGen.\qed
\smallni
{\it Remark.\/} As we mentioned the tensor $T=\Delta(z)^k\Delta(w)^kf(z,w)(dzdw)^{8k}$
can have poles along the 48 exceptional divisors. The results of Theorem \DegGen\
could be improved if one could find $f$ such that $T$ is holomorphic on the the whole
$\tilde\calB$. We did not succeed to find such $f$ and we think that they don't exist.
But we could not prove this.
\neupara{A two-fold covering of the box variety}%
We consider the subgroup $\Gamma'[4]$ of index two of $\Gamma[4]$ which is defined by
$$a+b+ c\equiv1\;\mod\; 8.$$
\proclaim
{Lemma}
{The group $\Gamma'[4]/\Gamma[8]$ is isomorphic to
$\gz/2\gz\times\gz/2\gz$. It acts freely on $\overline{\hz/\Gamma[8]}$.}
ActF%
\finishproclaim
{\it Proof.\/} Let $b$ by a point of the extended upper half plane 
$\hz^*=\hz\cup\qz\cup\{\infty\}$. Assume that $M\in \Gamma'[4]-\Gamma[8]$ is  a matrix
that fixes $a$ mod $\Gamma[8]$. Then there exists an element $A\in\Gamma[8]$ such that
$M(a)=A(a)$. The matrix $N=A^{-1}M$ then fixes $a$. This matrix is also contained
in $\Gamma'[4]$ and not in $\Gamma[8]$. Modulo 8 it is one of the following three
$$\pmatrix{1&4\cr4&1},\quad\pmatrix{5&4\cr0&5},\quad\pmatrix{5&0\cr4&5}.$$
Since it has a fixed point, the absolute value of its trace
is bounded by 2. We treat the three cases separately. In the first case we
have
$$N=\pmatrix{1+8\alpha&4+8\beta\cr4+8\gamma&1+8\delta}.$$
The condition for the trace implies $\delta=-\alpha$. The determinant is 1. 
But the condition $\delta=-\alpha$ implies that the determinant is
$1-16$ mod 32. This is a contraction.
\smallskip
In the second case we have
$$N=\pmatrix{5+8\alpha&4+8\beta\cr8\gamma&5+8\delta}.$$
The condition for the trace now gives $\delta=-\alpha-1$. Now the determinant would be
congruent 4 mod 8 which is not possible. The same argument works in the third case.\qed
\smallskip
We consider the (non-singular) manifold
$$X:=\overline{\hz/\Gamma[8]}\times \overline{\hz/\Gamma[8]}/\Gamma'[4]$$
where $\Gamma'[4]$ acts diagonally. The inclusion $\Gamma'[4]\hookrightarrow\Gamma[4]$
gives a two fold covering $X\to\calB$
of the box variety. Locally around the 48 singularities of $\calB$ this looks
like the covering $\cz^2\to\cz^2/\pm$. One can desingularize the node at 0 be first blowing up
the origin in the covering $\cz^2$. The involution $(z,w)\mapsto(-z,-w)$ 
lifts to this resolution and the quotient is smooth. The same can be done globally.
We blow up $X$ at the 48 inverse images of the nodes of $\calB$.
This gives a manifold $\tilde X$. The Galois involution of $X$ over $\calB$ lifts to $\tilde X$
and the quotient $\tilde \calB$ is just the blow up of $\calB$ at the nodes.
So we have a commutative diagram
$$\matrix{\tilde X&\lo&\tilde\calB\cr\downarrow&&\downarrow\cr
X&\lo&\calB}.$$
The map $\tilde X\to \tilde\calB$ is ramified along the exceptional divisors
(48 lines). The existence of this covering and its uniqueness have been treated in the
paper [Be]. 
We call 
$$X=\overline{\hz/\Gamma[8]}\times \overline{\hz/\Gamma[8]}/\Gamma'[4]$$
the Beauville manifold.
The existence of $X$ is equivalent to the fact that
the exceptional divisor  is divisible by two in the Picard group
$\Pic(\tilde\calB)$. The uniqueness follows from the fact that this Picard group is
torsion free [ST]. We refer to [Be] for more interesting properties of the surface $X$.
Some of them can be easily derived form the modular picture.
\smallskip

%
%
\smallskip
Since $X$ has the product of two curves of genus $>1$ as unramified covering, the universal
covering of $X$ is the product $\ez\times\ez$ of two unit disks. 
Let $C$ be the Riemann sphere or a torus $\cz/L$ and let $\alpha:C\to X$ be a holomorphic
map. This map lifts to the universal coverings. Since every holomorphic map
$\cz\to\ez$ is constant, we obtain that $\alpha$ is constant. Hence $X$ contains no
rational or elliptic curve. This argument applies also to the regular locus
$\calB_{\hbox{\sevenrm reg}}$ of the box variety, since its universal covering
is the complement of a discrete subset in $\ez\times\ez$. Hence we obtain the following
result.
\proclaim
{Remark}
{Every rational or elliptic curve in the boy variety contains at least one node.}
RaElB%
\finishproclaim
\neupara{Relation to a Kummer variety}%
In this section we consider the $\qz$-structure of the box variety. 
It is the associated projective variety of the
algebra
$$B=\qz[W_1,W_2,W_3,Z_1,Z_2,Z_3,C]$$
(with the defining relations of the box variety).
We consider the involution $\sigma(Z_3)=-Z_3$ of the ring $B$.
It induces an involution of the box variety. The invariant ring is
$$B^\sigma=\qz[Z_1,Z_2,C,W_1,W_2,W_3]$$
with defining relations
$$\eqalign{
W_1^2+W_3^2&=Z_2^2,\cr
W_2^2+W_3^2&=Z_1^2,\cr
W_1^2+W_2^2+W_3^2&=C^2.\cr}$$
The associated projective variety is the quotient of the box variety
by $\sigma$.
This is a also a modular variety, since in the picture of Theorem \MIs\ transformation
$\sigma$ is induced by the modular substitution
$$(z,w)\loma (Tz,w).$$
This gives the following result.
\proclaim
{Lemma}
{The variety $\calB/\sigma$ is defined over $\qz$. Over $\cz$ it agrees
with the modular variety which
belongs to the subgroup of $\SL(2,\gz)\times\SL(2,\gz)$ that is
generated by $\Gamma[8]\times\Gamma[8]$ and the elements
$$(T,E),\quad (E,T),\quad (T',T'),\quad (R,R).$$}
QuoM%
\finishproclaim
We considered already in Proposition \KleinR\
the group $\Gamma_1[8]\cap\Gamma[4]$ and explained the structure of the ring of modular forms.
We defined 4 generators $a,b,c,d$ with defining relations
$$a^2=c^2+d^2,\quad
b^2=c^2-d^2.$$
Since these relations are defined over $\qz$, we can can consider this 
algebra \hbox{over $\qz$}
$$A(\Gamma_1[8]\cap\Gamma[4])= C\otimes_\qz\cz,\quad C:=\qz[a,b,c,d].$$
and obtain an elliptic curve $E$ over $\qz$.
One can compute its normal form \hbox{over $\qz$}:
$$y^2=x^3-x.$$
Its projective form is $y^2z=x^3-xz^2$. An explicit isomorphism is given by
$$x=a-b,\quad y=2d,\quad z=2c - a - b.$$
We consider the automorphisms
$$\tau(a,b,c,d)=(a,-b,c,d),\quad \rho(a,b,c,d)=(a,b,-c,-d)$$
of the algebra.
\smallskip
They correspond to the modular transformations
$T'$ and $R$ whereas $T$ acts as identity. The two transformations generate
a group $H\cong\gz/2\gz\times\gz/2\gz$.
\proclaim
{Lemma}
{The transformation $\rho$ is an involution without fixed point of the
elliptic curve $E$. Hence it is a translation by a two-torsion point.
The transformation $\tau$ is an involution with the fixed point $[\sqrt2,0,1,1]$.
Hence the following is true. If one considers $E$ as an elliptic curve over
$\qz(\sqrt2)$ with origin $[\sqrt2,0,1,1]$ then $\tau$ corresponds to the negation
$x\mapsto -x$.}
FixN%
\finishproclaim
%
\smallskip
We want to consider the product of two copies of this curve.
This is the projective variety associated with the graded algebra
$$C_2=\qz[a\otimes a, a\otimes b,\dots, d\otimes d].$$
In the modular picture we have to identify
$$
a\otimes a=\vartheta_{00}(z)\vartheta_{00}(w),\dots,
d\otimes d=\vartheta_{10}(2z)\vartheta_{10}(2w).$$
The group $H$ acts diagonally on $C_2$. The fixed ring is just $B^\sigma$.
Hence we get the following result.
\proclaim
{Proposition}
{There is a biholomorphic map
$$(E\times E)/H\Isom\calB/\sigma,$$
defined over the field of Gauss numbers.}
IsGa%
\finishproclaim
The variety $(E\times E)/H$ can be understood as follows.
We first take the quotient by the translation $\varrho$. This gives an
abelian variety over $\qz$.
$$X=(E\times E)/\varrho.$$
Then we take the quotient by $\tau$. If we extend the base field $\qz$ by $\sqrt2$,
and take $[\sqrt2,0,1,1]$ diagonally embedded as origin then $\tau$ corresponds to the
negation and
$$(E\times E)/H=X/\pm$$
appears as a Kummer variety. Hence over the field $\qz(\imag,\sqrt2)$ of eighth roots
of unity, the variety $\calB/\sigma$ can be identified with a Kummer variety.
\neupara{A moduli problem}%
We denote by $(\Sch/S)$ the category of schemes over a base scheme $S$. For $S=\Spec(A)$
we write $(\Sch/A)$.
We fix  an algebraic number field $K$. We also fix an embedding $K\hookrightarrow\cz$.
Let $E$ be an elliptic curve over a scheme $S\in(\Sch/K)$
and $N$ a natural number. We denote by
$E[N]$ the kernel of multiplication by $N$ from $E$ to  $E$. This is a group scheme over
$S$. If $T$ as a scheme over $S$ then the $T$-valued points are
$$E[N](T)=\kernel(E(T)\buildrel \cdot N\over\lo E(T)).$$
To each finite group $G$ there is associated a group scheme which we denote by the
same letter and which is defined by $G(S)=G$ for connected $S$. We call it the constant
group scheme associated to $G$.
It may be that $E[N](S)$ is isomorphic $(\gz/N\gz\times \gz/N\gz)(S)$.
A level
$N$ structure on $E$ then means the choice of an isomorphism $(\gz/N\gz\times \gz/N\gz)(S)\Isom E[N](S)$.
It extends to an ismorphism of  $E(N)$
to the constant group scheme associated to  $\gz/N\gz\times \gz/N\gz$.
\smallskip
We denote by $G_m$ the group scheme ``multiplicative group'' over $K$.
It is defined by
$G_m(S)=\calO(S)^*$ (multiplicative group). The kernel of powering by $N$ is
the group scheme $\mu_N$. Hence
$$\mu_N(S)=\{a\in \calO(S);\ a^N=1\}.$$
For an elliptic curve $E$ over $S\in(\Sch/K)$ there exists the Weil-pairing.
It associates to each $S$-scheme $T$ an alternating map
$$E[N](T)\times E[N](T)\lo \mu_N(T).$$
When
$K$ contains the cyclotomic field of $N$-th roots of unity, then $\mu_N$ is the constant
group scheme associated to the abstract group
$$\mu_N=\{\zeta\in\cz;\ \zeta^N=1\}\cong\gz/N\gz.$$
We also can consider the symplectic pairing
$$e_N:\gz/N\gz\times\gz/N\gz\lo\mu_N,\quad e((a_1,a_2),(b_1,b_2))=
e^{2\pii(a_1b_2-a_2b_1)/N}.$$
From now on we assume that $K$ contains the $N$th roots of unity. Since we consider
a fixed embedding of $K$ into $\cz$, we can identify $\mu_N$ and $\mu_N(K)$.
So it makes sense to consider
level $N$-structures which preserve the symplectic pairings
and we will use the notion ``Level $N$-structure'' from now on in this
restricted sense.
\smallskip
We want to formulate a moduli problem. For this we introduce the category \Ell\
(compare [Ka]).
Its objects are pairs of elliptic curves $E\to S$, $F\to S$ over a variable scheme
$S\in\Sch/K$. Morphisms are pairs of cartesian squares
$$\matrix{E_1&\lo &E\cr\downarrow &&\downarrow\cr S_1&\lo& S}\qquad
\matrix{F_1&\lo &F\cr\downarrow &&\downarrow\cr S_1&\lo& S}$$
This category is a gruppoid, i.e.~all morphims are isomorphisms.
A moduli problem on $\Ell$ is a functor $\calP$ from $\Ell$ into the category of sets.
It is called representable if there exists a universal $\bfE\to\calM(\calP)$
in the sense that there is a functorial ismorphism
$$\calP(E/S)=\Hom_{\Ell}(E/S,\bfE/\calM(\calP)).$$
We recall that the moduli problem $\calP$ is called {\it relatively representable\/} if
for every $(E,F)/S$ in $\Ell$ the functor on $(\Sch/S)$ that associates
to a scheme $T/S$ the set $\calP((E_T,F_T)/T)$ is representable by some scheme
$\calP_{E/S}$.
\smallskip
Every representable moduli problem is also relatively representable. Let $\calP'$ and
$\calP''$ two moduli problems. Then one can define the simultaneous moduli problem
$\calP=(\calP',\calP'')$ in the obvious way. Assume that $\calP'$ is representatable
and $\calP''$ is relatively representable. Then it is easy to see that $\calP$
is representable.
\smallskip
We now consider two natural numbers $N\vert N'$.
We recall that we assume that $K$ containes the $N$th roots of unity.
We define a moduli problem $\calP(N,N')$.
It associates to $(E,F)/S$ the set of pairs $(\alpha,\beta)$
$$\alpha:\gz/N\gz\times\gz/N\gz\Isom E[N],\quad \beta: E[N']\Isom F[N'].$$
Here $\alpha$ is a level $N$-structure (preserving the Weil-pairing) and
$\beta$ is an isomorphism of $S$-groups, also preserving Weil-pairings.
We notice that $\alpha$ and $\beta$ induce a level $N$-structure of $F$.
\proclaim
{Theorem}
{Assume $N\ge 3$. The moduli problem $\calP(N,N')$ is representable.}
IsRep%
\finishproclaim
{\it Proof.\/} The moduli problem $\calP=\calP(N,N')$ can be considered as a
simultaneous problem $(\calP',\calP'')$. Here $\calP'$ refers to the level 4 moduli problem
and $\calP''$ associates to $(E,F)/S$ the set of isomorphisms $E[N']\to F[N']$
(preserving the Weil pairing). It is known that $\calP'$ is representable.
Hence it remains to show that $\calP''$ is locally representable.
For this we have
to fix a pair of  elliptic curves $E\to S$, $F\to S$ over a scheme
$S\in(\Sch/K)$. Then we have to consider the functor
on the category of schemes over  $S$ that associates to an scheme $T$ over $S$
the set of isomorphisms $E_T[N']\to F_T[N']$.
It is well-known that this functor is representable and also the subfunctor preserving
the Weil pairing is representable.\qed
\smallskip
We denote the moduli space of the functor $\calP(N,N')$ by
$\calM(N,N')$. This is an affine algebraic variety over $\qz(\imag)$.
Over $\cz$ it is
$\hz\times\hz/\Delta(N,N')$.
\smallskip
We are interested in the group $\Delta(4,8)$. Theorem \MIs\ says that
the surface $\hz\times\hz/\Delta(4,8)$ 
is embedded as an open part of the
box variety. We call this the {\it finite part of the box variety.\/}
We have now two $K$-structures, one coming form the defining equations
of the box variety and the other coming form the moduli problem $\calM(4,8)$.
A variant of the $q$-expansion principle shows that both are the same. Hence
we get the following result.
\proclaim
{Theorem}
{The variety $\calM(4,8)$ is isomorphic as variety over $\qz(\imag)$ to the
finite part of the box variety (considered as variety over $\qz(\imag)$).}
FinB%
\finishproclaim
This gives us a description of $\qz(\imag)$-rational points of the box variety
in terms of elliptic curves.
\proclaim
{Theorem}
{The $\qz(\imag)$-valued points of the finite part of the box variety are in one-to-one
correspondence to isomorphy classes of pairs
of elliptic curves $E,F$ over $\qz(\imag)$, equipped with a
level 4 structure of $E$ and a compatible isomorphism of group schemes $E[8]\to F[8]$ 
which preserves the Weil pairing.}
HSig%
\finishproclaim

\vskip1.5cm\noindent
{\paragratit References}%
\bigskip
\item{[Be]} Beauville, A.: {\it A tale of two surfaces,}, {http://arxiv.org/abs/1303.1910}
\medskip

\item{[Br]} Brieger, S.: {\it Modulformen zur Stufe 8},
Diplomarbeit, Universit\"at Heidelberg (2010)
\medskip

\item{[Fr]} Freitag, E.: {\it Eine Bemerkung zur Theorie der Hilbertschen Modulmannigfaltigkeiten
hoher Stufe,\/} Math. Z. {\bf 171}, 27--38 (1980)
\medskip

\item{[Ka]} Kani, E.: {\it Mazur's question and modular diagonal quotient surfaces,\/}
preprint (2011)
\medskip
\item{[Kl]} Klein, M.: {\it Ein Ring von Modulformen zur Stufe 8},
Diplomarbeit, Universit\"at Heidelberg (2011)
\medskip
\item{[vL]} van Luijk, R.: {\it On perfect cuboids,\/}
Undergraduate thesis, University Utrecht (2000)
\medskip
\item{[ST]}Stoll, M.,  Testa, D.: {\it The surface parametrizing cuboids}\hfill\break
 preprint   arXiv: 1009.0388
\bye